\documentclass[reqno]{amsart}
\usepackage{amsmath,amssymb}
\usepackage{bm}
\usepackage[dvipdfmx]{graphicx}
\usepackage[caption=false]{subfig}
\usepackage{verbatim}
\usepackage{wrapfig}
\usepackage{ascmac}
\usepackage{multirow}
\usepackage{amsfonts}
\usepackage{enumerate}
\usepackage{mathrsfs}
\usepackage{amsthm}
\usepackage{mathtools}
\mathtoolsset{showonlyrefs=true}
\usepackage[dvipdfmx]{xcolor}
\theoremstyle{plain}
\newtheorem{thm}{Theorem}
\newtheorem{prop}{Proposition}[section]
\newtheorem{lem}[prop]{Lemma}

\theoremstyle{remark}
\newtheorem{rem}[prop]{Remark}
\makeatletter
\renewenvironment{proof}[1][\proofname]{\par\pushQED{\qed}
  \normalfont
  \topsep6\p@\@plus6\p@ \trivlist
  \item[\hskip\labelsep{\bfseries #1}\@addpunct{\bfseries}]\ignorespaces
}{%
  \popQED\endtrivlist\@endpefalse
}
\renewcommand{\proofname}{Proof.}
\makeatother

\newcommand{\abs}[1]{\left\lvert#1\right\rvert}
\newcommand{\nabs}[1]{\lvert#1\rvert}
\newcommand{\norm}[1]{\left\lVert#1\right\rVert}
\newcommand{\R}{\mathbb{R}}
\newcommand{\E}{\mathcal{E}}
\newcommand{\T}{\mathcal{T}}

\def\coloneqq{\mathrel{\mathop:}=}%
\newcommand{\pderiv}[3][]{\frac{\partial^{#1}#2}{\partial #3^{#1}}}

\newcommand{\dist}{\operatorname{dist}}
\definecolor{cyan20}{cmyk}{.2,0,0,0}
\definecolor{MyBlue}{RGB}{15,82,188}

\makeatletter
    
    \@addtoreset{equation}{section}
  \makeatother

\begin{document}
\title[Nitsche's method for a Robin BVP]{
Nitsche's method for a Robin boundary value problem in a smooth domain
}
\author[Y. Chiba]{Yuki Chiba}
\address{Graduate School of Mathematical Sciences, The University of
Tokyo, Komaba 3-8-1, Meguro, Tokyo 153-8914, Japan}
\email{ychiba@ms.u-tokyo.ac.jp} 
\author[N. Saito]{Norikazu Saito}
\address{Graduate School of Mathematical Sciences, The University of
Tokyo, Komaba 3-8-1, Meguro, Tokyo 153-8914, Japan}
\email{norikazu@g.ecc.u-tokyo.ac.jp}
\urladdr{http://www.infsup.jp/saito/index-e.html}
\date{\today}

\subjclass[2000]{
Primary~
65N15, 	
Secondary~
65N30  	
}
\keywords{
finite element method, 
Nitsche's method
}
\maketitle

\begin{abstract}
We prove several optimal-order error estimates for a finite-element method applied to an inhomogeneous Robin boundary value problem (BVP) for the Poisson equation defined in a smooth bounded domain in $\mathbb{R}^n$, $n=2,3$. The boundary condition is weakly imposed using Nitsche's method. The Robin BVP is interpreted as the classical penalty method with the penalty parameter $\varepsilon$. The optimal choice of the mesh size $h$ relative to $\varepsilon$ is a non-trivial issue. This paper carefully examines the dependence of $\varepsilon$ on error estimates. 
Our error estimates require no unessential regularity assumptions on the solution. Numerical examples are also reported to confirm our results. 
\end{abstract}

\section{Introduction}
\label{sec:intro}

Nitsche's method \cite{nit71} is well-known as a powerful method for imposing the Dirichlet boundary condition (DBC) in the finite element method (FEM).  DBC is usually imposed by specifying the function values themselves at boundary nodal points. In contrast, Nitsche's method is based on the method of ``weak imposition''  of DBC using penalty parameter. Actually, this strategy is useful for resolving the issue of spurious oscillations for non-stationary Navier--Stokes and convection--diffusion equations as was mentioned in Bazilevs et al. \cite{bh07,bmch07}. 

In recent years, demand for computing complex boundary conditions has been increasing. Boundary conditions involving the Laplace--Beltrami operator, such as a dynamic boundary condition and a generalized Robin boundary condition play important roles in application to the reduced fluid--structure interaction model and Cahn--Hilliard equation (see, e.g., \cite{MR2243323}, \cite{MR2385883} and \cite{MR2629535}). 
Nitsche's method may be an effective approach to address these boundary conditions, and therefore, is worthy of a thorough investigation. 

When numerically solving PDEs in a smooth domain, 
we often utilize polyhedral approximations of the domain. 
Generally, a facile approximation of the problem may result in a wrong numerical solution; the so-called Babu\v{s}ka's paradox in \cite[\S 5]{MR0170133} is a remarkable example. Therefore, 
investigating not only the error caused by discretizations but also that caused by domain approximations is important. 
For the standard FEM, approximating domains is a common problem, and analysis of the energy norm is well-developed thus far. 
Only recently, the optimal order $W^{1,\infty}$ and $L^\infty$ stability and error estimates were established; refer to \cite{MR4071825} for detail.  

Consequently, we evaluate Nitsche's method for PDEs 
in a smooth bounded domain. 
We study a finite-element method (FEM) applied to an inhomogeneous Robin boundary value problem (BVP) for the Poisson equation defined in a smooth bounded domain $\Omega$ in $\mathbb{R}^n$, $n=2,3$. The boundary condition on the boundary $\Gamma=\partial\Omega$ is weakly imposed using Nitsche's method. We then derive several optimal-order error estimates under reasonable regularity assumptions on the solution. 
Specifically, we consider 
\begin{subequations}
  \label{eq:poi}
\begin{alignat}{2}
 -\Delta u &=f &{ }\qquad&\mbox{in}\  \Omega \label{eq:poi1}\\
 \displaystyle\pderiv{u}{\nu} + \frac{1}{\varepsilon}u&=\displaystyle \frac{1}{\varepsilon}u_0 + g &&  \mbox{on}\ \Gamma.  \label{eq:poi2}
\end{alignat}
Therein, we suppose that $f\in L^2(\Omega)$,  $u_0\in H^{2}(\Omega)$, and $g\in H^{1}(\Omega)$ are given. 
Moreover, $\partial/\partial\nu$ denotes the differentiation along the outward unit normal vector $\nu$ to $\Gamma$ and $\varepsilon>0$ is a constant. 
If $\varepsilon=0$, we consider \eqref{eq:poi1} with the Dirichlet boundary condition
\begin{equation}
 u=u_0 \qquad \mbox{on}\ \Gamma   \label{eq:poi3}
\end{equation}
\end{subequations}
instead of \eqref{eq:poi2}. 

The case of a polyhedral domain with $\varepsilon>0$ has already been addressed in Juntunen and Stenberg \cite{MR2501054}. We are motivated by \cite{MR2501054} and this paper is a generalization of \cite{MR2501054} to a smooth domain. 
We study simultaneously the case $\varepsilon=0$, that is the case of DBC.

If we are concerned with the Dirichlet BVP \eqref{eq:poi1} and \eqref{eq:poi3}, the Robin BVP \eqref{eq:poi1} and \eqref{eq:poi2} with $g=0$ implies the classical penalty method with the penalty parameter $\varepsilon\to 0$. 
(The $\varepsilon$ is interpreted as the penalty parameter in the classical penalty method. On the other hand, Nitsche's method is introduced using the penalty parameter, which we will write as $\gamma$. The readers have to care not to confuse.)
FEM for this method is well studied so far; we refer to \cite{bab73, be86, MR968097,kin74} for example. In particular, Barret and Elliott \cite{be86} presented the error estimate in a smooth domain as
\begin{equation}
 \label{eq:be86}
\|\tilde{u}-u_h^{\varepsilon}\|_{H^1(\Omega_h)}
\le C
[(h+\varepsilon^{-\frac12}h^2+\varepsilon^{\frac12}) \|u\|_{H^4(\Omega)}
+\varepsilon^{-\frac12}h^2\|\tilde{u}_0\|_{H^3(\Omega_h)}
].
\end{equation}
Therein, $h$ denotes the granularity parameter of the triangulation and 
$\Omega_h$ denotes a polyhedral approximation of $\Omega$ satisfying \eqref{eq:assumption}. (See also Remark \ref{rem:assumption}.) 
The continuous $\mathcal{P}^1$ finite element solution is represented by $u_h^{\varepsilon}$. The definition of function spaces and their norms are described in the end of this section. Moreover, $\tilde{u}$ and $\tilde{u}_0$ are suitable extensions of $u$ and $u_0$, respectively. The precise definition will be mentioned in the next section. 
The estimate \eqref{eq:be86} gives the optimal-order estimate for the $H^1$ norm by setting $\varepsilon=Ch^2$. However, we need a surplus regularity $u\in H^4(\Omega)$. Barret and Elliott \cite{MR968097} later studied the iso-parametric FEM for a similar problem and obtained similar results as ours. 
However, regularity assumptions slightly vary from ours.    

This paper carefully examines the dependence of $\varepsilon$ on error estimates. As a matter of fact, how to choose $h$ relative to $\varepsilon$ is a non-trivial issue for smooth domain cases. A suitable regularity assumption is another non-trivial issue as we recalled for the standard FEM above. This point was not discussed in Juntunen and Stenberg \cite{MR2501054}, because these issues do not appear in polyhedral cases. In fact, we succeed in deriving (see Theorem \ref{thm1})
\begin{equation}
 \norm{\tilde u - u_h}_{h,\star}
\le 
C h \left[\left(1+\frac{h}{\varepsilon}\right)\norm{u}_{2,\Omega}+\norm{f}_\Omega +\frac{h}{\varepsilon}\norm{u_0}_{2,\Omega}+\norm{g}_{1,\Omega}\right] ,
\label{eq:errorh1z}
\end{equation} 
where $\norm{\cdot}_{h,\star}$ denotes the DG $H^1$ norm defined as \eqref{eq:norm_n}. Consequently, we deduce the optimal-order estimate for the DG $H^1$ norm by setting $h=C\varepsilon$ under no further assumptions the smoothness of the solution and data. It makes this possible by applying some estimates reported in \cite{MR4071825}. On the other hand, we assume a surplus regularity $u\in W^{3,q}(\Omega)$, $q>n$, for deducing the optimal-order estimate for the $L^2$ norm (see Theorem \ref{thm2}). In our opinion, this is an essential requirement; see Remark \ref{rem:reg} and Section \ref{sec:numerical_example}.

This paper comprises six sections. 
In Section \ref{sec:model}, the continuous $\mathcal{P}^1$ finite element approximation using Nitsche's method and the main error estimates (Theorems \ref{thm1} and \ref{thm2}) are described. 
After having presented some preliminary results in Section \ref{sec:preliminaries}, we prove Theorems \ref{thm1} and \ref{thm2} in 
Sections \ref{sec:pr_thmh1} and \ref{sec:pr_thml2}, respectively.  
Finally, numerical examples are also reported to confirm our results in Section \ref{sec:numerical_example}.

\subsection*{Notation} 
We list the notations used in this paper. 
We follow the standard notation of, for example, \cite{MR2424078} for function spaces and their norms. In particular, for $1\le p \le \infty$ and a positive integer $j$, we use the standard Lebesgue space $L^{p}(\mathcal{O})$ and Sobolev space $W^{j,p}(\mathcal{O})$. Hereinafter, $\mathcal{O}$ denotes a bounded domain in $\R^n$. The inner product and norm of  $L^2(\mathcal{O})$ are denoted, respectively, by $(\cdot,\cdot)_{\mathcal{O}}$ and $\|\cdot\|_{\mathcal{O}}$. 
The norm of $W^{j,p}(\mathcal{O})$ is denoted by $\|\cdot\|_{W^{j,p}(\mathcal{O})}$. As usual, we set $H^j(\mathcal{O})=W^{j,2}(\mathcal{O})$, and the semi-norm and norm of $H^{j}(\mathcal{O})$ are denoted by, respectively,  
\[
    \abs{v}_{i,\mathcal{O}} = \left(\sum_{\abs{\alpha} = i}\norm{\pderiv[\alpha]{v}{x}}_{\mathcal{O}}^2\right)^{{1}/{2}},\quad 
    \norm{v}_{j,\mathcal{O}} = \left(\sum_{i = 0}^{j} \abs{v}_{i,\mathcal{O}}^2\right)^{{1}/{2}}.
\]
For $S\subset \partial \mathcal{O}$, we define $L^{p}(S)$ 
using a surface measure $dS$ in a common approach.
The inner product and norm of $L^2(S)$ is denoted by, respectively,  
$\langle\cdot,\cdot\rangle_{S}$ and $\|\cdot\|_{S}$.
Moreover, $\mathcal{P}^r(\mathcal{O})$ denotes the set of all polynomials in $\mathcal{O}$ of degree $\le r$.

\section{Nitsche's method and the main results}
\label{sec:model}

We recall that $\Omega$ is a bounded domain  in $\mathbb{R}^n$, $n=2,3$. 
Throughout this paper, we assume that the boundary $\Gamma$ of $\Omega$ is a $C^{l}$ boundary, where $l$ is an integer $\ge 2$. 

From the general theory of elliptic PDEs, we know that the unique solution $u$ of \eqref{eq:poi} belongs to $H^2(\Omega)$ and satisfies 
$\|u\|_{2,\Omega}\le C_\varepsilon (\|f\|_{\Omega}+\|u_0\|_{2,\Omega}+\|g\|_{1,\Omega})$, where $C_\varepsilon$ denotes a positive constant depending only on $\Omega$ and $\varepsilon$. 
The pure Neumann problem is out of our interest. Therefore, we assume 
\begin{equation}
\label{eq:ep0}
0\le \varepsilon \le \varepsilon_0
\end{equation}  
for a suitably large $\varepsilon_0>0$. 

Let $\{\mathcal{T}_h\}_{h}$ be a regular family of
triangulations $\mathcal{T}_h$ of a polyhedral domain $\Omega_h\subset \R^n$ in the sense of \cite{Cia78}. That is,  
\begin{enumerate}
\item $\mathcal{T}_h$ is a set of closed $n$-simplices (elements) $K$, and 
\[
 \Omega_h=\operatorname{Int}\left(\bigcup_{K\in\mathcal{T}_h}K\right);
\] 
\item The granularity parameter $h$ is defined as $h =  \displaystyle \max_{K \in \T_h}h_K$, where $h_K$ denotes the diameter of $K$;   
\item Any two elements of $\mathcal{T}_h$ meet only in entire common
      faces or sides or in vertices;
\item There exists a positive constant $c_{\textup{reg}}$ satisfying  
$h_K \le c_{\textup{reg}} \rho_K$ for all $K\in\mathcal{T}_h\in \{\mathcal{T}_h\}_{h}$, where $\rho_K$ denotes the diameter of the inscribed ball of $K$. 
\end{enumerate}
We then introduce the boundary mesh $\mathcal{E}_h$ inherited from $\mathcal{T}_h$ by 
\[
 \E_h\coloneqq\{E\subset \partial\Omega_h \colon  \mbox{$E$ is an $(n-1)$-face of some $K\in\mathcal{T}_h$}\},
\]
and the boundary $\Gamma_h\coloneqq\partial \Omega_h$ is expressed as $\Gamma_h = \bigcup_{E \in \E_h} E$. We assume that $\Gamma_h$ is an approximate surface/polygon of $\Gamma$ in the sense that 
\begin{equation}
\label{eq:assumption} 
\mbox{every vertex of $E \in \E_h$ lies on $\Gamma$.}
\end{equation} 

We use the continuous $\mathcal{P}^1$ finite element space
\begin{equation}
  V_h  \coloneqq 
\{ \chi \in C(\overline \Omega) 
\colon \chi|_K \in \mathcal{P}^1(K)\ (\forall K\in \T_h)\}. 
\label{eq:fespace_n}
\end{equation}

Below we fix a sufficiently smooth domain $\tilde \Omega\subset \R^n$ such that 
\begin{equation}
 \overline{\Omega\cup\Omega_h}\subset \tilde\Omega.
\end{equation}
Since $\Gamma$ is of class $C^l$, $l\ge 2$, the domain $\Omega$ admits a strong $l$-extension operator $P$. That is, $P$ is a linear operator of $W^{l,r}(\Omega)\to W^{l,r}(\tilde \Omega)$ for any $0\le k\le l$ and $1\le r< \infty$, and it satisfies 
\begin{equation}
\label{eq:extension}
 (Pv)|_\Omega=v\mbox{ in }\Omega,\quad 
\|Pv\|_{W^{k,r}(\tilde{\Omega})}\le C_{l,r,\Omega}\|v\|_{W^{k,r}(\Omega)},
\end{equation}
where $C_{l,r,\Omega}$ denotes a positive constant depending only on $l$, $r$ and $\Omega$; see \cite[Theorem 5.22]{MR2424078} for example. 
Using this, we write 
\begin{equation}
\tilde{f}\coloneqq Pf,\quad 
\tilde{u}_0\coloneqq Pu_0,\quad 
\tilde{g}\coloneqq Pg.
 \label{eq:ext}
\end{equation}

Following \cite{MR2501054}, we set
\begin{subequations} 
  \begin{align}
 a_h(w,v) &\coloneqq  (\nabla w,\nabla v)_{\Omega_h} + b_h(w,v),\\
 b_h(w,v)&\coloneqq  \sum_{E \in \E_h}\left[-\frac{\gamma h_E}{\varepsilon+\gamma h_E}\left(\langle\frac{\partial w}{\partial \nu_h},v\rangle_{E}+\langle w,\pderiv{v}{\nu_h}\rangle_{E}\right) \right.\nonumber\\
       &\hspace{7em}\left.+ \frac{1}{\varepsilon+\gamma h_E}\langle w,v\rangle_E -\frac{\varepsilon\gamma h_E}{\varepsilon+\gamma h_E}\langle\pderiv{w}{\nu_h},\pderiv{v}{\nu_h}\rangle_E
        \right],
     \end{align}
and 
\begin{multline}
      l_h(v) \coloneqq  (\tilde f, v)_{\Omega_h} +
       \sum_{E \in \E_h}\left[\{\frac{1}{\varepsilon+\gamma h_E}\langle \tilde u_{0},v\rangle_{E} -\frac{\gamma h_E}{\varepsilon+\gamma h_E}\langle \tilde u_0, \pderiv{v}{\nu_h}\rangle_E \right. \\
       \left. + \frac{\varepsilon}{\varepsilon+\gamma h_E}\langle \tilde g,v\rangle_{E} -\frac{\varepsilon\gamma h_E}{\varepsilon+\gamma h_E}\langle \tilde g, \pderiv{v}{\nu_h}\rangle_E
       \right],
      \end{multline}
\end{subequations}
where $\gamma>0$ is a penalty parameter, $h_E$ the diameter of $E$, and $\nu_h$ the outer unit normal vector to $\Gamma_h$.  

Nitsche's method for \eqref{eq:poi} is stated as follows: 
      \begin{equation}
        u_h \in V_h, \qquad  
        a_h(u_h,\chi) = l_h(\chi) \quad (\forall \chi \in V_h).
        \label{eq:nitsche}
      \end{equation}

We use the following norms that depend on $\varepsilon$ and $h_E$: 
      \begin{align}
        \norm{v}_{h}^2 &\coloneqq \norm{\nabla v}_{\Omega_h}^2 + \sum_{E \in \E_h} \frac{1}{\varepsilon+\gamma h_E}\norm{v}_{E}^2,\\
        \norm{v}_{h,\star}^2 &\coloneqq \norm{v}_h^2 + \sum_{E \in \E_h} h_E\norm{\pderiv{v}{\nu_h}}_{E}^2, \label{eq:norm_n}
      \end{align} 
We recall that $\norm{\cdot}_{h}$ and $\norm{\cdot}_{h,\star}$ are equivalent on $V_h$ uniformly in $h$. That is, there exists a positive constant $C_0$ independent of $h$ such that 
\begin{equation}
\label{eq:equiv}
\|v\|_h\le \|v\|_{h,\star}\le C\|v\|_{h}\qquad (v\in V_h).
\end{equation}
Here and hereinafter, $C$ denotes a generic positive constant which is independent of $h$ and $\varepsilon$. The value of $C$ may be different at each occurrence. 
The inequalities \eqref{eq:equiv} follow from the well-known inequality
\begin{equation}
\label{eq:3.1}
 \sum_{E\in\mathcal{E}_h}h_E\left\|\frac{\partial v}{\partial \nu_h}\right\|_E^2\le C\|\nabla v\|_{K}^2\qquad (v\in V_h,~K\in\mathcal{T}_h,~E\subset K). 
\end{equation} 
In fact, \eqref{eq:3.1} is a readily obtainable consequence of the standard trace inequality, 
\begin{equation}
\|v\|_{E}^2 \le Ch_E^{-1}\left(\|v\|_K^2+h_K^2\|\nabla v\|_K^2\right) \qquad  (E\in\mathcal{E}_h,E\subset\partial K, v\in H^1(K)).
\label{eq:trace0}
\end{equation}

Nitsche's method \eqref{eq:nitsche} admits a unique solution in view of 
the following basic result; see \cite[Theorem 3.2]{MR2501054}. 

\begin{lem}
We have 
\begin{subequations} 
  \begin{equation}
    a_h(w,v) \le C \norm{w}_{h,\star}\norm{v}_{h,\star}\qquad (\forall w,\,v \in H^2(\Omega_h) + V_h). \label{eq:continuity}
  \end{equation}
Moreover, there exists a positive constant $\gamma_0$ which is independent of $h$ and $\varepsilon$ such that we have for $0<\gamma\le \gamma_0$,  
  \begin{equation}
    a_h(\chi,\chi) \ge C \norm{\chi}_h^2 \qquad (\forall \chi \in V_h). \label{eq:coercivity}
  \end{equation}
\end{subequations}
\end{lem}

Actually, $\gamma_0$ can be taken as any positive number strictly smaller than $1/C$, where $C$ denotes the constant appearing in \eqref{eq:3.1}; see \cite[Theorem 3.2]{MR2501054}. Below we always assume that 
\begin{equation}
 \label{eq:pen9}
0<\gamma\le \gamma_0.
\end{equation}
To deduce convergence results, we need an inverse assumption as
\begin{equation}
\label{eq:inv1}
 h \le C \min_{E\in\mathcal{E}_h}h_E.  
\end{equation}

We are now in a position to state our main result. 
We recall that $C$ denotes a positive constant which is independent of $\varepsilon$ and $h$. 

\begin{thm}[$H^1$ estimates]\label{thm1}
Suppose that $\Gamma$ is a $C^2$ boundary. 
Let $u \in H^2(\Omega)$ and $u_h \in V_h$ represent the solutions of 
\eqref{eq:poi} and \eqref{eq:nitsche}, respectively. 
Assume that 
\eqref{eq:assumption}, \eqref{eq:pen9} and 
\eqref{eq:inv1} are satisfied. 
Then, if $\varepsilon>0$, 
we have  
\begin{equation}
 \norm{\tilde u - u_h}_{h,\star}
\le 
C h \left[\left(1+\frac{h}{\varepsilon}\right)\norm{u}_{2,\Omega}+\norm{f}_\Omega +\frac{h}{\varepsilon}\norm{u_0}_{2,\Omega}+\norm{g}_{1,\Omega}\right] ,
\label{eq:errorh1}
\end{equation}
where $\tilde{u}=Pu$. If $\varepsilon=0$, 
we have
\begin{equation}
 \norm{\tilde u - u_h}_{h,\star}
\le 
C h (\norm{u}_{2,\Omega}+\norm{f}_\Omega + \norm{u_0}_{2,\Omega}) .
\label{eq:errorh1a}
\end{equation}
\end{thm}

\begin{thm}[$L^2$ estimates]\label{thm2}
Suppose that $\Gamma$ is a $C^3$ boundary. 
Let $u \in H^2(\Omega)$ and $u_h \in V_h$ represent the solutions of 
  \eqref{eq:poi} and \eqref{eq:nitsche}, respectively.
Assume that 
\eqref{eq:assumption}, \eqref{eq:pen9} and 
\eqref{eq:inv1} are satisfied. 
Then, if $\varepsilon>0$, $u\in W^{3,q}(\Omega)$ for some $q>n$, $u_0\in H^3(\Omega)$ and $g\in H^3(\Omega)$, we have  
  \begin{align}
    \norm{\tilde u - u_h}_{\Omega_h} &\le Ch^2 \Biggl[
      \norm{u}_{W^{3,q}(\Omega)} + \left(\frac{h}{\varepsilon} +\frac{h^2}{\varepsilon^2}\right)\norm{u}_{2,\Omega}  +\norm{u_0}_{3,\Omega} \nonumber\\
      &\hspace{5em}+ \left(\frac{h}{\varepsilon} +\frac{h^2}{\varepsilon^2}\right)\norm{u_0}_{2,\Omega}
      +\norm{g}_{3,\Omega} + \frac{h}{\varepsilon}\norm{g}_{1,\Omega}
      \Biggr] , \label{eq:errorl2}
  \end{align}
  where $\tilde{u}=Pu$. On the other hand, if $\varepsilon=0$ and $u\in H^3(\Omega)$,  we have
  \begin{equation}
   \norm{\tilde u - u_h}_{\Omega}
  \le Ch^2 (\norm{u}_{3,\Omega} + \norm{u_0}_{2,\Omega}) .
  \label{eq:errorl2a}
  \end{equation}
\end{thm}

\begin{rem}
\label{rem:reg}
Theorem \ref{thm1} reports that the optimal rate of convergence for the $H^1$ error is achieved under a reasonable (minimal) regularity assumption on $u$. On the other hand, we pose a somewhat surplus regularity $u\in W^{3,q}(\Omega)$, $q>n$, for deducing the optimal rate of convergence for the $L^2$ error.  
In our opinion, this is an essential requirement. Actually, a numerical example reported in Section \ref{sec:numerical_example} shows the second-order convergence may not take place if $u\in W^{3,q}(\Omega)$, $q>n$. 
\end{rem}

\begin{rem}
\label{rem:assumption}
We are assuming \eqref{eq:assumption} for $\Omega_h$. This can be replaced by 
\begin{equation}
 \operatorname{dist}(\Omega,\Omega_h)\le Ch^2
\end{equation}
with some obvious modification of proofs. 
\end{rem}

\section{Boundary-skin estimates}
\label{sec:preliminaries}

We collect some auxiliary results that will be used in the proof of the main results. 

Since $\Omega$ is a bounded $C^l$ domain, there exists a local coordinate system $\{U_r,y_r,\phi_r\}_{r=1}^M$ to ensure the following:
    \begin{enumerate}[1)]
        \item $\{U_r\}_{r=1}^M$ is an open covering of $\Gamma$.\label{enum_dom1}
        \item For any $1\le r\le M$, there exists a congruent transformation $A_r$ such that $y_r= (y_{r1},y'_{r}) = A_r(x)$, where $x$ is the original coordinate.\label{enum_dom2}
        \item For any $1\le r\le M$, $\phi_r$ is a $C^l$ function in $\Delta_r\coloneqq \{y_{r1} \in \R^{n-1} \colon \abs{y_{r1}}\le \alpha\} $ and $\Gamma \cap U_r$ is a graph of $\phi_r$ with respect to the coordinate $y_r$.\label{enum_dom3}
    \end{enumerate}
Assuming that $h$ is sufficiently small if necessary, our possible assumptions are as follows:  
\begin{enumerate}[1)]
\setcounter{enumi}{3}
 \item For any $1\le r\le M$, there exists a function $\phi_{rh}$ such that $\Gamma_h \cap U_r$ is a graph of $\phi_{rh}$ with respect to the coordinate $y_r$.
\end{enumerate}

In addition, we assume that $h_0$ is sufficiently small to ensure that for any $x \in \Gamma$ and $r=1,\ldots,M$, the open ball $B(x,h_0)$ with center $x$ and radius $h_0$ is contained in a neighborhood $U_r$. 
Let $d(x)$ be the signed distance function defined by
\[
  d(x) \coloneqq \begin{cases}
    -\dist(x,\Gamma) & x \in \Omega \\
    \dist(x,\Gamma) & x \in \mathbb{R}^n\backslash\Omega.
  \end{cases}
  \]
  We define $\Gamma(\delta) \coloneqq \{x \in \R^n \colon \abs{d(x)} < \delta\}$, which we call the \emph{boundary-skin} region. 
  Then, for a sufficiently small $\delta$, the orthogonal projection $\pi$ onto $\Gamma$ exists such that
  \begin{equation}
    x = \pi(x) + d(x) \nu(\pi(x)) \quad (x \in \Gamma(\delta),\,\pi(x) \in \Gamma). 
  \end{equation} 
 
Because $h$ is sufficiently small, $\pi$ is defined on $\Gamma_h\subset \Gamma(\delta)$ and for each $E\in \E_h$, and $\pi(E)$ comprises some local neighborhood $U_r$.
  In this case, $\pi|_{\Gamma_h}$ has the inverse operator $\pi^*(x)=x+t^*(x)\nu(x)$, and $\norm{t^*}_{L^{\infty}(\Gamma)}\le Ch^2$. Moreover, $\pi(\E_h)\coloneqq\{\pi(E) \colon E \in \E_h\}$ is a partition of $\Gamma$.

We assume that all these properties hold for any $h\le h_0$ by assuming that $h_0$ is sufficiently small if necessary. 


Now we can state the boundary-skin estimates. For the proof, refer to \cite[Theorems 8.1, 8.2, and 8.3 and Lemma 9.1]{MR3563279} and \cite[Lemma A.1]{MR4071825}. 

\begin{lem}[Boundary-skin estimates]\label{lem:bdskin}
Let $\delta \in [c^{-1}h^2,~ ch^2]$ with a positive constant $c$. 
We have 
\begin{subequations} 
\begin{equation}
     \norm{\nu_h-\nu\circ\pi}_{L^\infty(\Gamma_h)} \le Ch.  \label{eq:bs_5} 
\end{equation}
For $f \in W^{1,p}(\Gamma(\delta))$, we have
  \begin{align}
 \norm{f-f\circ\pi}_{L^p(\Gamma_h)} &\le C\delta^{1-1/p}\norm{f}_{W^{1,p}(\Gamma(\delta))} ,\label{eq:bs_2}\\
        \norm{f}_{L^p(\Gamma(\delta))} &\le C(\delta \norm{\nabla f}_{L^p(\Gamma(\delta))}+\delta^{1/p}\norm{f}_{L^p(\Gamma)}) .\label{eq:bs_3}
      \end{align}
Moreover, for $f \in W^{1,p}(\Omega_h)$, we have 
\begin{equation}
\norm{f}_{L^p(\Omega_h\setminus\Omega)} \le C(\delta \norm{\nabla f}_{L^p(\Omega_h\setminus\Omega)}+\delta^{1/p}\norm{f}_{L^p(\Gamma_h)}) .\label{eq:bs_4}\\
\end{equation}
\end{subequations}
\end{lem}

\begin{lem}
\label{la:skin1}
\begin{equation}
\label{eq:skin1a}
\|\chi\|_{\Omega_h\backslash\Omega}
\le Ch \|\chi\|_{h}\qquad (\chi\in V_h)
\end{equation}
\end{lem}

\begin{proof}
It is a direct consequence of \eqref{eq:bs_4} in view of \eqref{eq:ep0} and \eqref{eq:pen9}.  
\end{proof}

\begin{lem}
\label{la:skin2}
\begin{subequations}
 If $\varepsilon=0$, we have  
\begin{equation}
\label{eq:skin1c}
\norm{\tilde u-\tilde u_0}_{\Gamma_h}
\le 
Ch^2\left(\|u\|_{2,\Omega}+\|u_0\|_{2,\Omega}\right) .
\end{equation}
Moreover, if $\varepsilon>0$, we have 
\begin{equation}
\label{eq:skin1b}
\norm{\pderiv{\tilde u}{\nu_h} + \frac{\tilde u-\tilde u_0}{\varepsilon}-\tilde g}_{\Gamma_h}
\le 
Ch\left(\|u\|_{2,\Omega}+
\frac{h}{\varepsilon}\|u\|_{2,\Omega}
+\frac{h}{\varepsilon}\|u_0\|_{2,\Omega}+\|g\|_{1,\Omega}\right).
\end{equation}
\end{subequations}
\end{lem}

\begin{proof}
First, consider the case $\varepsilon=0$. 
Let $\delta$ be a constant satisfying $c^{-1}h^2\le \delta\le ch^2$ as in Lemma \ref{lem:bdskin}. 
Since $u\circ\pi=u\circ u_0$ on $\Gamma_h$, we have by \eqref{eq:bs_2}   
\begin{align*}
\norm{\tilde u-\tilde u_0}_{\Gamma_h} &\le 
\norm{\tilde u-\tilde u\circ \pi}_{\Gamma_h}+\norm{\tilde u_0\circ\pi-\tilde u_0}_{\Gamma_h} \\
&\le 
Ch\|\tilde{u}\|_{1,\Gamma(\delta)}+Ch\|\tilde{u}_0\|_{1,\Gamma(\delta)}\\
&\le 
Ch^2\|\tilde{u}\|_{\tilde{\Omega}}+Ch^2\|\tilde{u}_0\|_{\tilde\Omega},
\end{align*}
which implies \eqref{eq:skin1c}. Therein, we have used
$\|\tilde{u}\|_{1,\Gamma(\delta)}\le Ch\|\tilde{u}\|_{2,\Omega}$ and 
$\|\tilde{u}_0\|_{1,\Gamma(\delta)}\le Ch\|\tilde{u}_0\|_{2,\Omega}$, 
which are direct consequences of \eqref{eq:bs_3} and the trace theorem. 

We proceed to the case $\varepsilon>0$.  
Since $(\nabla \tilde{u}\circ\pi)\cdot (\nu\circ\pi)+\varepsilon^{-1}\tilde{u}\circ\pi=\tilde{g}\circ\pi +\varepsilon^{-1}\tilde{u}_0\circ\pi$ on $\Gamma_h$, we have  
\begin{align*}
\norm{\pderiv{\tilde u}{\nu_h} + \frac{\tilde u-\tilde u_0}{\varepsilon}-\tilde g}_{\Gamma_h}
&\le \|\nabla \tilde{u}\cdot \nu_h-(\nabla \tilde{u}\circ\pi)\cdot (\nu\circ\pi)\|_{\Gamma_h}+\norm{\tilde g-\tilde g\circ \pi}_{\Gamma_h}
\\
&\mbox{ }\quad +\frac{1}{\varepsilon}\norm{\tilde u-\tilde u\circ \pi}_{\Gamma_h}
+\frac{1}{\varepsilon}\norm{ \tilde{u}- \tilde{u}\circ \pi}_{\Gamma_h}\\
&=:\textup{J}_1+\textup{J}_2+\textup{J}_3+\textup{J}_4.
\end{align*}
As above, we estimate as
\[
 \textup{J}_2+\textup{J}_3+\textup{J}_4\le 
\frac{h^2}{\varepsilon}\|\tilde{u}\|_{2,\Omega}+
\frac{h^2}{\varepsilon}\|\tilde{u}_0\|_{2,\Omega}+\|g\|_{1,\Omega}.
\]
On the other hand, by \eqref{eq:bs_2} and \eqref{eq:bs_5}
\begin{align*}
 \textup{J}_1 &\le \|\nabla \tilde{u}\cdot \nu_h-\nabla \tilde{u}\cdot (n\circ\pi)\|_{\Gamma_h} +\|\nabla \tilde{u}\cdot (\nu\circ\pi)-
(\nabla \tilde{u}\circ\pi)\cdot (\nu\circ\pi)\|_{\Gamma_h}\\
&\le Ch\|\tilde{u}\|_{2,\Omega}.
\end{align*}
Summing up, we deduce \eqref{eq:skin1b}. 
\end{proof}

\section{Proof of Theorem \ref{thm1}}
\label{sec:pr_thmh1}

We start with a version of the Strang lemma. To state it, we set
\begin{equation}
 \label{eq:res}
r_h(v,\chi)=a_h(v,\chi)-l_h(\chi)
\end{equation}
for $v\in H^2(\tilde{\Omega})+V_h$ and $\chi\in V_h$. 

\begin{lem}
Under the same assumption of Theorem \ref{thm1}, we have
\begin{equation}
  \norm{\tilde u - u_h}_{h,\star} \le C 
\left[
\inf_{\xi\in V_h}\norm{\tilde u - \xi}_{h,\star}+\sup_{\chi \in V_h} \frac{\nabs{r_h(\tilde{u},\chi)}}{\norm{\chi}_{h}}
\right].
\label{eq:approx_h1}
\end{equation}
\end{lem}

\begin{proof}
 Letting $\xi \in V_h$ and $\chi = u_h-\xi$, we have
\begin{align*}
\norm{\chi}_h^2 &\le C_2 a_h(\chi,\chi) \\
  &= C[a_h(\tilde u-\xi,\chi)-a_h(\tilde u,\chi)+ l_h(\chi)]\\
  & \le C \norm{\tilde u -\xi}_{h,\star}\norm{\chi}_{h,\star} + \nabs{r_h(\tilde{u},\chi)},
\end{align*}
where \eqref{eq:coercivity}, \eqref{eq:continuity}, and \eqref{eq:equiv} are applied. This, together with the triangular inequality, implies \eqref{eq:approx_h1}
\end{proof}

The standard Lagrange interpolation operator of $C(\overline{\Omega})\to V_h$ is denoted by $\Pi_h$. It is well-known that 
\begin{equation}
\abs{w-\Pi_h w}_{m,K} \le C h_K^{2-m}|w|_{2,K} \quad (w \in H^2(\Omega_h), ~ K \in \T_h, ~ m = 0,\,1,\,2). 
\label{eq:interpolation}
\end{equation}
As a direct application of \eqref{eq:interpolation} and \eqref{eq:trace0}, we derive 
\begin{equation}
\|w-\Pi_hw\|_{h,\star}\le Ch|w|_{2,\Omega}\qquad (w\in H^2(\Omega)).
\label{eq:lag}
\end{equation}

We can state the following proof. 

\begin{proof}[Proof of Theorem \ref{thm1}.]
First consider the case $\varepsilon>0$. In view of \eqref{eq:approx_h1} and \eqref{eq:lag}, it suffices to prove 
\begin{equation}
 \sup_{\chi\in V_h}\frac{\abs{r_h(\tilde u,\chi)}}{\norm{\chi}_{h}}\le 
Ch\left(\norm{u}_{2,\Omega}+\frac{h}{\varepsilon}\norm{u}_{2,\Omega}+\|f\|_\Omega+\frac{h}{\varepsilon}\norm{u_0}_{2,\Omega}+\norm{g}_{1,\Omega}\right).
\label{eq:errorh10}
\end{equation}

Let $\chi\in V_h$ be arbitrary.   
Applying the integrate by parts, we have
\begin{align}
a_h(\tilde u,\chi) - l_h(\chi) &= (-\Delta \tilde u-\tilde{f},\chi)_{\Omega_h} \nonumber \\
&\mbox{ }\qquad + \sum_{E \in \E_h}\frac{\varepsilon}{\varepsilon + \gamma h_E}\langle \pderiv{\tilde u}{\nu_h} + \frac{\tilde u-\tilde u_0}{\varepsilon}-\tilde g,\chi - \gamma h_E \pderiv{\chi}{\nu_h} \rangle_E. 
 \nonumber \\ 
& =:\textup{I}_1+\textup{I}_2. \label{eq:pfthm1h11}
\end{align}
Since $-\Delta \tilde u = \tilde f$ in $\Omega$, we have by \eqref{eq:skin1a}
\begin{align*}	
\abs{\textup{I}_1} &=\abs{(-\Delta \tilde u-\tilde{f},\chi)_{\Omega_h\backslash\Omega}}\\
&\le \norm{\Delta \tilde u+ \tilde{f}}_{\Omega_h\setminus\Omega}\norm{\chi}_{\Omega_h\setminus\Omega}  
\le C h ( \norm{u}_{2,\Omega}+\norm{f}_{\Omega} )\norm{\chi}_{h}. 
\end{align*}

Using \eqref{eq:skin1b}, 
\begin{align}
|\textup{I}_2|  & \le C\norm{\pderiv{\tilde u}{\nu_h} + \frac{\tilde u-\tilde u_0}{\varepsilon}-\tilde g}_{L^2(\Gamma_h)}\norm{\chi}_{h,\star} \nonumber \\
& \le Ch\left(\norm{u}_{2,\Omega}+\frac{h}{\varepsilon}\norm{u}_{2,\Omega}+\frac{h}{\varepsilon}\norm{u_0}_{2,\Omega}+\norm{g}_{1,\Omega}\right)\norm{\chi}_{h}. 
\label{eq:pfthm1h1_3}
\end{align}
Summing up, we deduce \eqref{eq:errorh10}.  
 
We proceed to the case $\varepsilon=0$. It suffices to prove
\begin{equation}
 \sup_{\chi\in V_h}\frac{\abs{r_h(\tilde u,\chi)}}{\norm{\chi}_{h}}\le 
Ch\left(\norm{u}_{2,\Omega}+\|f\|_\Omega+\norm{u_0}_{2,\Omega}\right).
\label{eq:errorh1a0}
\end{equation} 
In this case, we have $r_h(\tilde u,\chi) =\textup{I}_1+\textup{I}_3$, where 
\[
 \textup{I}_3=\sum_{E \in \E_h}\frac{1}{\gamma h_E}\langle \tilde u-\tilde u_0,\chi - \gamma h_E \pderiv{\chi}{\nu_h} \rangle_E. 
\]
We apply \eqref{eq:inv1} and \eqref{eq:skin1c} to obtain
\[
 |\textup{I}_3| 
\le Ch_{\min}^{-1}\|\tilde u-\tilde u_0\|_{\Gamma_h}\|\chi\|_{h,\star} 
\le Ch(\|u\|_{2,\Omega}+\| u_0\|_{\Omega})\|\chi\|_{h},
\]
where $h_{\min}=\min_{E}h_E$. Therefore, \eqref{eq:errorh1a0} is proved. 
\end{proof}

\section{Proof of Theorem \ref{thm2}}
\label{sec:pr_thml2}

We use the Green operator $G:L^2(\Omega_h)\to H^2(\Omega)$ defined as follows. 
For $\varepsilon>0$, $z=G\eta\in H^2(\Omega)$ is the unique solution of  
\begin{subequations}
\label{eq:z}
\begin{gather}
-\Delta z  = 
\begin{cases}
 \eta & \mbox{in }\Omega\cap \Omega_h\\
0 & \mbox{in }\Omega\backslash \Omega_h,
\end{cases}
 \label{eq:zeq}\\ 
\pderiv{z}{\nu} + \frac{1}{\varepsilon}z =0 \quad\mbox{on } \Gamma, \label{eq:zbc}
\end{gather}
where $\eta\in L^2(\Omega_h)$. For $\varepsilon =0$, \eqref{eq:zbc} is replaced by 
\begin{equation}
 z=0\quad \mbox{on } \Gamma.
\label{eq:zbc0}
\end{equation} 
\end{subequations}
For any $\eta\in L^2(\Omega_h)$, $G\eta$ admits the a priori estimate
\begin{equation}
\label{eq:green}
 \|G\eta\|_{2,\Omega}\le C\|\eta\|_{\Omega_h};
\end{equation}
In particular, $C$ does not depend on $\varepsilon$. We omit the proof since it is outside the scope of this paper.
As a matter of fact, this can be verified by a standard method of difference quotient. For example, if tracing the proof of \cite[Theorem 3.3]{MR3296617} carefully, we find that the estimate \eqref{eq:green} holds true. Moreover, if $\Gamma$ is a $C^4$ boundary, the same proof of \cite[Lemma 4.1]{sai04} is also applicable.  

\begin{lem}
\label{lem:approx2}
Under the same assumption of Theorem \ref{thm2}, we have
\begin{multline}
\norm{\tilde u - u_h}_{L^2(\Omega_h)} 
\le C\Biggl[ \norm{\tilde u - u_h}_{L^2(\Omega_h\setminus\Omega)} + hc_1(h,\varepsilon)\norm{\tilde u - u_h}_{h,\star} \\
+ \sup_{\eta \in L^2(\Omega_h)}\frac{\norm{\tilde z-\Pi_h\tilde z}_{h,\star}\norm{\tilde u-u_h}_{h,\star}+\nabs{r_h(\tilde u,\Pi_h\tilde z)}}{\norm{z}_{2,\Omega}}\Biggr], 
\label{eq:approxl2}
\end{multline}
where $\tilde z = Pz\in H^2(\tilde{\Omega})$, $z=G{\eta}\in H^2(\Omega)$ and 
\[
 c_1(h,\varepsilon)=
\begin{cases}
 1+h/\varepsilon & (\varepsilon>0),\\
 1 & (\varepsilon=0).
\end{cases}
\] 
\end{lem}

\begin{proof}
First, suppose that $\varepsilon>0$. In the similar way as the derivation of \eqref{eq:pfthm1h11}, we deduce 
\begin{equation}
\label{eq:consistency}
(\chi,-\Delta v)_{\Omega_h}  
=
a_h(\chi,v)-\sum_{E \in \E_h}
\frac{\varepsilon}{\varepsilon+\gamma h_E}
\langle \chi-\gamma h_E\pderiv{\chi}{\nu_h},\pderiv{v}{\nu_h}+\frac{1}{\varepsilon}v \rangle_E
\end{equation}
for $\chi \in H^2(\Omega_h)+V_h$ and $v \in H^2(\Omega_h)$.

    
Let $\eta\in L^2(\Omega_h)$ and $\tilde{z}=Pz=P(G\eta)$. 
We use the same symbol $\eta$ to express the zero extension of $\eta$ into $\Omega\backslash\Omega_h$. 
Substituting \eqref{eq:consistency} for $\chi = \tilde u -u_h$ and $v = \tilde z$, we obtain
    \begin{align*}
 (\tilde u-u_h,\eta)_{\Omega_h} 
&= (\tilde u - u_h,-\Delta \tilde z)_{\Omega_h} + (\tilde u - u_h,\eta + \Delta \tilde z)_{\Omega\setminus\Omega_h} \nonumber \\
      &= a_h(\tilde u -u_h,\tilde z) + (\tilde u - u_h,\eta + \Delta \tilde z)_{\Omega\setminus\Omega_h} \nonumber\\
      &\hspace{5em}- \sum_{E \in \E_h}\frac{\varepsilon}{\varepsilon+\gamma h_E}\langle \tilde u -u_h-\gamma h_E\pderiv{(\tilde u -u_h)}{\nu_h},\pderiv{\tilde z}{\nu_h}+\frac{1}{\varepsilon}\tilde z \rangle_E\nonumber\\
      &= a_h(\tilde u -u_h,\tilde z-\Pi_{\#}\tilde z)+r_h(\tilde u,\Pi_h\tilde z) + (\tilde u - u_h,\eta + \Delta \tilde z)_{\Omega\setminus\Omega_h} \nonumber\\
      &\hspace{5em}- \sum_{E \in \E_h}\frac{\varepsilon}{\varepsilon+\gamma h_E}\langle \tilde u -u_h-\gamma h_E\pderiv{(\tilde u -u_h)}{\nu_h},\pderiv{\tilde z}{\nu_h}+\frac{1}{\varepsilon}\tilde z \rangle_E,\nonumber \\
&=:\textup{J}_1+\textup{J}_2+\textup{J}_3+\textup{J}_4.\label{eq:pf_approxl2_1}
    \end{align*}
Thanks to  \eqref{eq:green}, an estimation for $\textup{J}_3$ is readily; 
\[
 \textup{J}_3 \le C
 \norm{\tilde u - u_h}_{\Omega\setminus\Omega_h}\norm{\eta + \Delta \tilde z}_{\Omega\setminus\Omega_h} \le C\norm{\tilde u - u_h}_{\Omega\setminus\Omega_h}\|\eta\|_{\Omega_h}.
\]
We apply Lemma \ref{la:skin2} with $u_0=0$ and $g=0$ and derive
    \begin{align*}
      \textup{J}_4 &\le C \norm{\tilde u -u_h}_{h,\star}\norm{\pderiv{\tilde z}{\nu_h}+\frac{1}{\varepsilon}\tilde z}_{L^2(\Gamma_h)}\nonumber \\
      &\le C h\left(\norm{z}_{2,\Omega}+\frac{h}{\varepsilon}\|z\|_{2,\Omega}\right)\norm{\tilde u -u_h}_{h,\star}. \label{eq:pf_approxl2_2}
    \end{align*}
Summing up, we deduce
    \begin{align*}
      \norm{\tilde u - u_h}_{L^2(\Omega_h)} & = \sup_{\eta \in L^2(\Omega_h)} \frac{(\tilde u-u_h,\eta)_{\Omega_h}}{\norm{\eta}_{L^2(\Omega_h)}} \nonumber \\
      & \le C\norm{\tilde u - u_h}_{L^2(\Omega\setminus\Omega_h)}+Ch\left(1+\frac{1}{\varepsilon}\right)\norm{\tilde u -u_h}_{h,\star} \nonumber \\
      & \hspace{3em} C \sup_{\eta \in L^2(\Omega_h)} \frac{|a_h(\tilde u -u_h,\tilde z-\Pi_{\#}\tilde z)|+|r_h(\tilde u,\Pi_h\tilde z)|}{\norm{\eta}_{L^2(\Omega_h)}}.
    \end{align*}
Finally, using \eqref{eq:continuity} and \eqref{eq:green}, we deduce \eqref{eq:approxl2} for $\varepsilon>0$. 

If $\varepsilon=0$, $\textup{J}_4$ is replaced by 
\[
 \textup{J}_4'=
- \sum_{E \in \E_h}\frac{1}{\gamma h_E}
\langle \tilde u -u_h-\gamma h_E 
\pderiv{(\tilde u -u_h)}{\nu_h},\tilde z \rangle_E.
\]
We apply Lemma \ref{la:skin2} with $u_0=0$ and get
\[
 |\textup{J}_4'|\le Ch_E^{-1}\norm{\tilde u -u_h}_{h,\star}\cdot Ch^2\norm{{z}}_{2,\Omega}.
\]
Therefore, \eqref{eq:approxl2} holds even for $\varepsilon=0$.  
\end{proof}

We finally state the following proof. 

\begin{proof}[Proof of Theorem \ref{thm2}.]
We define following bilinear and linear forms: 
\begin{subequations} 
\begin{align}
  a(w,v) &= (\nabla w,\nabla v)_{\Omega} + b(w,v)\\
   b(w,v)&= \sum_{E \in \E_h}\biggl\{-\frac{\gamma h_E}{\varepsilon+\gamma h_E}\Bigl(\langle\pderiv{w}{\nu},v\rangle_{\pi(E)}+\langle w,\pderiv{v}{\nu}\rangle_{\pi(E)}\Bigr),\nonumber\\
    &\hspace{7em}+ \frac{1}{\varepsilon+\gamma h_E}\langle w,v\rangle_{\pi(E)} -\frac{\varepsilon\gamma h_E}{\varepsilon+\gamma h_E}\langle\pderiv{w}{\nu},\pderiv{v}{\nu}\rangle_{\pi(E)}
     \biggr\},\\
l(v) &= (\tilde f, v)_{\Omega} +
       \sum_{E \in \E_h}\biggl\{\frac{1}{\varepsilon+\gamma h_E}\langle \tilde u_{0},v\rangle_{\pi(E)} -\frac{\gamma h_E}{\varepsilon+\gamma h_E}\langle \tilde u_0, \pderiv{v}{\nu}\rangle_{\pi(E)}\nonumber\\
       &\hspace{10em}+ \frac{\varepsilon}{\varepsilon+\gamma h_E}\langle \tilde g,v\rangle_{\pi(E)} -\frac{\varepsilon\gamma h_E}{\varepsilon+\gamma h_E}\langle \tilde g, \pderiv{v}{\nu}\rangle_{\pi(E)}
       \biggr\}
  \end{align}
\end{subequations}
  Then, we obtain
  \[
    a(u,v) = l(v) \quad (\forall v \in H^s(\Omega)).
    \]
    First consider the case $\varepsilon>0$.
    Since $\tilde u$ and $\tilde z$ are continuous in $\Omega_h$, we have 
    \begin{align}
      a_h(\tilde u,\tilde z)-l_h(\tilde z) &=a_h(\tilde u,\tilde z) - a(\tilde u,\tilde z) + l(\tilde z) -l_h(\tilde z) \nonumber \\
      &= \left(\int_{\Omega_h\setminus\Omega}(\nabla \tilde u\cdot \nabla \tilde z - \tilde f\tilde z)\,dx-\int_{\Omega\setminus\Omega_h}(\nabla \tilde u \cdot \nabla \tilde z - \tilde f \tilde z)\,dx\right) \nonumber \\
      &{ }\quad +\sum_{E\in\E_h}\frac{1}{\varepsilon+\gamma h_E} \biggl[\langle -\gamma h_E \pderiv{\tilde u}{\nu_h}+\tilde u - \tilde u_0 - \varepsilon \tilde g,\tilde z \rangle_E \nonumber\\
      &\hspace{10em}-\langle -\gamma h_E \pderiv{u}{\nu}+u - u_0 - \varepsilon g, z \rangle_{\pi(E)}\biggr] \nonumber\\
      &{ }\quad -\sum_{E\in\E_h}\frac{\varepsilon\gamma h_E}{\varepsilon+\gamma h_E}\langle  \pderiv{\tilde u}{\nu_h}+\frac{\tilde u  - \tilde u_0}{\varepsilon} - \tilde g,\pderiv{\tilde z}{\nu_h}\rangle_E\nonumber \\
      &=: \textup{I}_1+\textup{I}_2-\textup{I}_3. \label{eq:prthm1l2_1}
    \end{align}
    Using \eqref{eq:bs_3}, we obtain 
    \begin{align}
      \abs{\textup{I}_1}&\le Ch^2\left(\norm{\nabla\tilde u \cdot \nabla\tilde z}_{W^{1,1}(\widetilde \Omega)} + \norm{\tilde f \tilde z}_{W^{1,1}(\widetilde \Omega)} \right) \nonumber \\
      &\le Ch^2\norm{u}_{3,\Omega}\norm{z}_{2,\Omega} \label{eq:prthm1l2_2}
    \end{align}

Given that
\[
  \langle w,v \rangle_E -\langle w,v \rangle_{\pi(E)}  = \int_E(wv-(w\circ\pi)(v\circ\pi))\,d\gamma_h + \int_E (w\circ\pi)(v\circ\pi)\,d\gamma_h - \int_{\pi(E)} wv\,d\gamma,
\]
we have
\begin{align*}
  \textup{I}_2 &= \sum_{E\in\E_h}\frac{1}{\varepsilon+\gamma h_E}\Bigl[\langle-\gamma h_E \pderiv{\tilde u}{\nu_h}+\tilde u-\tilde u_0-\varepsilon\tilde g,\tilde z\rangle_{E}
  -\langle -\gamma h_E \pderiv{u}{\nu}+u-u_0-\varepsilon g,z\rangle_{\pi(E)}\Bigr]  \\
  &= \sum_{E\in\E_h}\frac{1}{\varepsilon+\gamma h_E}\Bigl[\langle -\gamma h_E \nabla \tilde u \cdot (\nu_h-\nu\circ\pi), \tilde z \rangle_{E} \\
  &\qquad+ \langle-\gamma h_E \nabla \tilde u \cdot \nu\circ\pi+\tilde u-\tilde u_0-\varepsilon\tilde g,\tilde z\rangle_{E}
  -\langle -\gamma h_E \nabla u\cdot \nu+u-u_0-\varepsilon g,z\rangle_{\pi(E)}\Bigr] \\
  &= \sum_{E\in\E_h}\frac{-\gamma h_E }{\varepsilon+\gamma h_E} \langle  \nabla \tilde u \cdot (\nu_h-\nu\circ\pi), \tilde z \rangle_{E}
   + \sum_{E\in\E_h}\frac{-\gamma h_E }{\varepsilon+\gamma h_E} \Bigl[\langle \nabla \tilde u \cdot \nu\circ\pi,\tilde z\rangle_{E} - \langle \nabla u\cdot \nu,z\rangle_{\pi(E)}\Bigr] \\
   &\qquad+ \sum_{E\in\E_h}\frac{\varepsilon}{\varepsilon+\gamma h_E} \Bigl[\langle \frac{\tilde u - \tilde u_0}{\varepsilon},\tilde z\rangle_{E} -\langle \frac{u-u_0}{\varepsilon},z\rangle_{\pi(E)}\Bigr]
   - \sum_{E\in\E_h}\frac{\varepsilon}{\varepsilon+\gamma h_E} \Bigl[\langle \tilde g,\tilde z\rangle_{E} -\langle g,z\rangle_{\pi(E)}\Bigr]\\
   &=: \textup{A}_1 + \textup{A}_2 + \textup{A}_3 - \textup{A}_4.
\end{align*}
Using \eqref{eq:bs_5}, we obtain
\begin{align*}
  \abs{\textup{A}_1} &\le C \frac{h}{\varepsilon}\sum_{E \in \E_h} \norm{ \nabla \tilde u \cdot (\nu_h-\nu\circ\pi) \tilde z }_{L^1(E)} \\
  & \le C \frac{h^2}{\varepsilon} \norm{\tilde u}_{L^2(\Gamma_h)}\norm{\tilde z}_{L^2(\Gamma_h)} \le C\frac{h^2}{\varepsilon} \norm{u}_{2,\Omega}\norm{\tilde z}_{L^2(\Gamma_h)}.
\end{align*}
Using \eqref{eq:zbc} and \eqref{eq:bs_2}, we have
\begin{align*}
  \norm{\tilde z}_{L^2(\Gamma_h)} & \le \norm{\tilde z-z\circ \pi}_{L^2(\Gamma_h)}+\norm{z\circ \pi}_{L^2(\Gamma_h)} \\
  &\le Ch\norm{z}_{1,\Omega} + C \norm{z}_{L^2(\Gamma)} \\
  & \le Ch\norm{z}_{1,\Omega} + C\norm{\varepsilon \nabla z \cdot \nu}_{L^2(\Gamma)} \\
  & \le C(\varepsilon + h) \norm{z}_{2,\Omega}.
\end{align*}
Hence, we obtain
\[\abs{\textup{A}_1} \le Ch^2(1+\frac{h}{\varepsilon})\norm{u}_{2,\Omega}\norm{z}_{2,\Omega}.\]
We rearrange $\textup{A}_2$ as
\begin{align*}
  \textup{A}_2 &= \sum_{E \in \E_h}\frac{-\gamma h_E }{\varepsilon+\gamma h_E} \int_E (\nabla \tilde u \cdot \nu) \circ \pi) (\tilde z- z\circ\pi)  d\gamma_h\\
   &\qquad+\sum_{E \in \E_h}\frac{-\gamma h_E }{\varepsilon+\gamma h_E} \int_E ((\nabla \tilde u -\nabla u\circ \pi) \cdot \nu \circ \pi) (z\circ \pi)~ d\gamma_h\\
    &\qquad+ \sum_{E \in \E_h}\frac{-\gamma h_E }{\varepsilon+\gamma h_E} \Bigl[\int_E ((\nabla u)\circ \pi \cdot \nu \circ \pi) (z\circ\pi) ~ d\gamma_h -  \int_{\pi(E)} (\nabla u \cdot \nu) z ~ d\gamma\Bigr] \\
    &= \textup{A}_{21}+\textup{A}_{22} + \textup{A}_{23}.
\end{align*}
Using boundary skin estimates, we get
\begin{gather*}
  \abs{\textup{A}_{21}} \le C h \norm{\nabla \tilde u \cdot \nu \circ \pi }_{L^2(\Gamma_h)}\norm{\tilde z}_{H^1(\Gamma(\delta))} \le Ch^2 \norm{u}_{2,\Omega}\norm{z}_{2,\Omega}, \\
  \abs{\textup{A}_{22}} \le C h^2 \norm{\nabla \tilde u}_{W^{1,\infty}(\Gamma(\delta))}\norm{z\circ \pi}_{L^1(\Gamma_h)} \le Ch^2 \norm{u}_{W^{3,q}(\Omega)}\norm{z}_{1,\Omega}, \\
  \abs{\textup{A}_{23}} \le C h^2 \norm{\nabla u \cdot \nu z}_{L^1(\Gamma)} \le C h^2 \norm{u}_{2,\Omega}\norm{z}_{1,\Omega}.
\end{gather*}
Therein, we used the Sobolev inequality and \eqref{eq:extension} to derive
\[
\norm{\tilde u}_{W^{2,\infty}(\tilde{\Omega})}
\le 
C \norm{\tilde u}_{W^{3,q}(\tilde{\Omega})}
\le 
C \norm{u}_{W^{3,q}(\Omega)}.
\]
(To apply \eqref{eq:extension} we are assuming that $\Gamma$ is a $C^3$ boundary.)

Therefore, we deduce 
\[\abs{\textup{A}_2} \le C h^2 \norm{u}_{W^{3,q}(\Omega)}\norm{z}_{2,\Omega}.\]
Similarly, we have
\[\abs{\textup{A}_4} \le C h^2 \norm{g}_{3,\Omega}\norm{z}_{2,\Omega}.\]
$\textup{A}_3$ is rewritten as
\begin{align*}
  \textup{A}_3 &= \sum_{E \in \E_h}\frac{\varepsilon}{\varepsilon+\gamma h_E} \int_E \Bigl[ \frac{\tilde u - \tilde u_0}{\varepsilon} (\tilde z- z\circ\pi) \Bigr] d\gamma_h\\
   &\qquad+\sum_{E \in \E_h}\frac{\varepsilon}{\varepsilon+\gamma h_E} \int_E \Bigl[(\frac{\tilde u - \tilde u_0}{\varepsilon}-\frac{u\circ \pi - u_0\circ \pi}{\varepsilon})z\circ \pi\Bigr]d\gamma_h\\
    &\qquad+ \sum_{E \in \E_h}\frac{\varepsilon}{\varepsilon+\gamma h_E} \Bigl[\int_E\frac{u\circ \pi - u_0\circ \pi}{\varepsilon} z\circ\pi d\gamma_h -  \int_{\pi(E)} \frac{u- u_0}{\varepsilon}z d\gamma\Bigr] \\
    &= \textup{A}_{31}+\textup{A}_{32} + \textup{A}_{33}.
\end{align*}
Using boundary skin estimates and \eqref{eq:zbc}, we can perform estimations as 
\begin{align*}
  \abs{\textup{A}_{31}} &\le C \frac{h}{\varepsilon} \norm{\tilde u - \tilde u_0}_{L^2(\Gamma_h)}\norm{\tilde z}_{H^1(\Gamma(\delta))} \\
& \le C\frac{h^2}{\varepsilon}\norm{\tilde u - \tilde u_0}_{L^2(\Gamma_h)}(h \norm{\tilde z}_{H^2(\Gamma(\delta))} + \norm{z}_{H^1(\Gamma)})\\
  &\le Ch^2(1+\frac{h}{\varepsilon}) (\norm{u}_{1,\Omega}+\norm{u_0}_{1,\Omega})\norm{z}_{2,\Omega},
\end{align*}
and 
\begin{align*}
  \abs{\textup{A}_{32}} \le C \frac{h^2}{\varepsilon} \norm{\tilde u - \tilde u_0}_{W^{1,\infty}(\Gamma(\delta))}\norm{z\circ \pi}_{L^1(\Gamma_h)} \le Ch^2(\norm{u}_{3,\Omega}+\norm{u_0}_{3,\Omega})\norm{z}_{2,\Omega}.
\end{align*}
Moreover, 
\begin{align*}
  \abs{\textup{A}_{33}} \le C \frac{h^2}{\varepsilon} \norm{(\tilde u - \tilde u_0) z}_{L^1(\Gamma)} \le C h^2 (\norm{u}_{1,\Omega}+\norm{u_0}_{1,\Omega})\norm{z}_{2,\Omega}.
\end{align*}
So, we get
\begin{align*}
  \abs{\textup{A}_3} \le C h^2\left(\norm{u}_{3,\Omega}+\frac{h}{\varepsilon}\norm{u}_{1,\Omega}+\norm{u_0}_{3,\Omega}+\frac{h}{\varepsilon}\norm{u_0}_{1,\Omega}\right)\norm{z}_{2,\Omega}.
\end{align*}
Summing up, we deduce
\begin{align}
  \abs{\textup{I}_2} \le Ch^2\left(\norm{u}_{W^{3,q}(\Omega)}+\frac{h}{\varepsilon}\norm{u}_{2,\Omega}+\norm{u_0}_{3,\Omega}+\frac{h}{\varepsilon}\norm{u_0}_{1,\Omega}+\norm{g}_{3,\Omega}\right)\norm{z}_{2,\Omega}.
\end{align}
We apply \eqref{eq:skin1b}, we obtain
\begin{align*}
  \abs{\textup{I}_3} \le Ch^2 \left(\norm{u}_{2,\Omega}+\frac{h}{\varepsilon}\norm{u}_{2,\Omega}+\frac{h}{\varepsilon}\norm{u_0}_{2,\Omega}+\norm{g}_{1,\Omega}\right)\norm{z}_{2,\Omega}.
\end{align*}
    Consequently, 
    \begin{multline*}
      \nabs{a_h(\tilde u,\tilde z)-l_h(\tilde z)} \\
\le Ch^2\left(\norm{u}_{W^{3,q}(\Omega)}+\frac{h}{\varepsilon}\norm{u}_{2,\Omega}+\norm{u_0}_{3,\Omega}+\frac{h}{\varepsilon}\norm{u_0}_{2,\Omega}+\norm{g}_{3,\Omega}\right)\norm{z}_{2,\Omega}.
    \end{multline*}

    Using the same way as the proof of \eqref{eq:errorh10}, we derive
    \[
     \nabs{a_h(\tilde u,w)-l_h(w)} \le Ch\left(\norm{u}_{2,\Omega}+\frac{h}{\varepsilon}\norm{u}_{2,\Omega}+\frac{h}{\varepsilon}\norm{u_0}_{2,\Omega}+\norm{g}_{1,\Omega}\right)\norm{w}_{h,\star}
    \]
    for $w \in H^2(\tilde \Omega) + V_h$.
    By substituting this for $w = \tilde z-\Pi_h\tilde z$, we have
    \begin{align*}
      \nabs{r_h(\tilde u,\Pi_h\tilde z)} &\le \nabs{a_h(\tilde u,\tilde z-\Pi_h\tilde z)-l_h(\tilde z-\Pi_h\tilde z)} +  \nabs{a_h(\tilde u,\tilde z)-l_h(\tilde z)} \nonumber \\
      &\le C  h^2\left(\norm{u}_{W^{3,q}(\Omega)}+\frac{h}{\varepsilon}\norm{u}_{2,\Omega}+\norm{u_0}_{3,\Omega}+\frac{h}{\varepsilon}\norm{u_0}_{2,\Omega}+\norm{g}_{3,\Omega}\right)\norm{z}_{2,\Omega}.
    \end{align*}
    Finally, we have
    \[
      \norm{\tilde u -u_h}_{\Omega_h\setminus\Omega} \le Ch\norm{\tilde u -u_h}_{h,\star}
    \]
    and we obtain the estimate \eqref{eq:errorl2}.

    We proceed to the case $\varepsilon = 0 $. We replace $\textup{I}_2$ by
    \begin{align*}
      \textup{I}_2' 
  &= -\sum_{E\in\E_h}\langle  \nabla \tilde u \cdot (\nu_h-\nu\circ\pi), \tilde z \rangle_{E}
   -\sum_{E\in\E_h}\Bigl[\langle \nabla \tilde u \cdot \nu\circ\pi,\tilde z\rangle_{E} - \langle \nabla u\cdot \nu,z\rangle_{\pi(E)}\Bigr] \\
   &\qquad+ \sum_{E\in\E_h}\frac{1}{\gamma h_E} \Bigl[\langle \tilde u - \tilde u_0,\tilde z\rangle_{E} -\langle u-u_0,z\rangle_{\pi(E)}\Bigr]\\
   &= \textup{A}_1' + \textup{A}_2 + \textup{A}_3',
    \end{align*}
    and $\textup{I}_3$ by
    \begin{align*}
      \textup{I}_3' 
  &= \sum_{E\in\E_h}\langle \tilde u  - \tilde u_0,\pderiv{\tilde z}{\nu_h}\rangle_E.
    \end{align*}
    Using \eqref{eq:bs_5} and \eqref{eq:bs_2}, we have
    \begin{align*}
      \nabs{\textup{A}_1'} &\le Ch \norm{\nabla \tilde u}_{L^2(\Gamma_h)}\norm{\tilde z - z \circ \pi}_{L^2(\Gamma_h)} \\
      & \le Ch^2 \norm{u}_{2,\Omega}\norm{z}_{2,\Omega}.
    \end{align*}
    Using \eqref{eq:skin1c} and \eqref{eq:bs_2}, we obtain
    \begin{align*}
      \nabs{\textup{A}_3'} &\le C h_{\min}\norm{\tilde u - \tilde u_0}_{L^2(\Gamma_h)}\norm{\tilde z - z \circ \pi}_{L^2(\Gamma_h)} \\
      & \le Ch^2 (\norm{u}_{2,\Omega}+\norm{u_0}_{2,\Omega})\norm{z}_{2,\Omega}.
    \end{align*}
    Since $\textup{A}_2 = \textup{A}_{21}$, we get
    \begin{align*}
      \nabs{\textup{A}_2} = Ch^2(\norm{u}_{2,\Omega}+\norm{u_0}_{2,\Omega})\norm{z}_{2,\Omega}.
    \end{align*}
    Using \eqref{eq:skin1c}, we have
    \begin{align*}
      \nabs{\textup{I}_3'} &\le \norm{\tilde u - \tilde u_0}_{L^2(\Gamma_h)}\norm{\pderiv{\tilde z}{\nu_h}}_{L^2(\Gamma_h)} \\
      & \le Ch^2 (\norm{u}_{2,\Omega}+\norm{u_0}_{2,\Omega})\norm{z}_{2,\Omega}.
    \end{align*}
    Hence, we obtain
    \begin{align}
      \nabs{r_h(\tilde u,\Pi_h\tilde z)} &\le \nabs{a_h(\tilde u,\tilde z-\Pi_h\tilde z)-l_h(\tilde z-\Pi_h\tilde z)} +  \nabs{a_h(\tilde u,\tilde z)-l_h(\tilde z)} \nonumber \\
      &\le C  h^2\left(\norm{u}_{3,\Omega}+\norm{u_0}_{2,\Omega}+\norm{g}_{2,\Omega}\right)\norm{z}_{2,\Omega}.
    \end{align}
    Therefore, \eqref{eq:errorl2a} is proved.
\end{proof}

\section{Numerical examples}
\label{sec:numerical_example}
In this section, we present some numerical results to verify the validity of our error estimates.
We consider the Poisson problem \eqref{eq:poi} in a disk $\Omega = \{x \in \R^2 \colon \abs{x} < 1\}$.

\begin{figure}[bt]
  \centering
  \subfloat[][$\varepsilon = 1$]{\includegraphics[scale = 0.55]{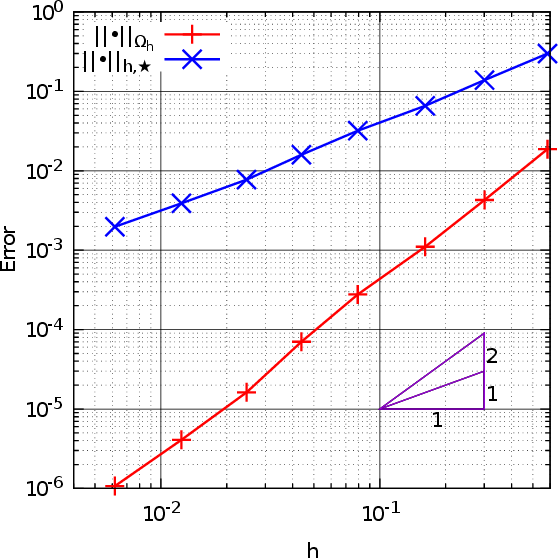}\label{fig:n_nsym_p1}}
     \subfloat[][$\varepsilon = 0$]{\includegraphics[scale = 0.55]{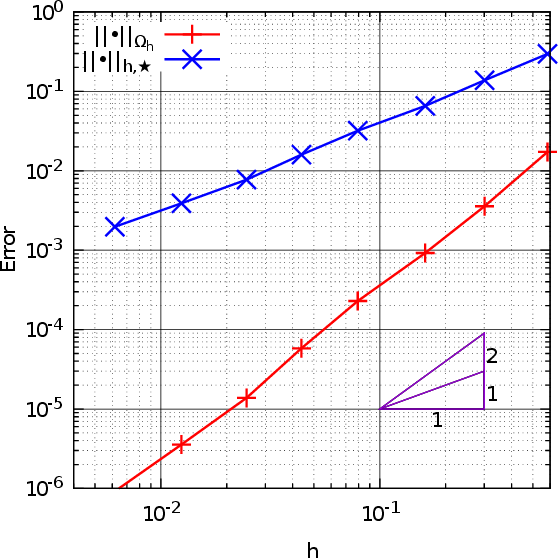}\label{fig:n_nsym_p1_dir}}
  \caption{$H^1$ errors and $L^2$ errors for the exact solution $u=\sin(x_1)\sin(x_2)$. }
  \label{fig:nsym_p1_err}
\end{figure}

First, we confirm the validity of the estimates in Theorems \ref{thm1} and \ref{thm2}. 
We set $f$, $g$ and $u_0$ so that a function $u(x_1,x_2) = \sin(x_1)\sin(x_2)$ solves \eqref{eq:poi}. Let $u_h\in V_h$ be the solution of \eqref{eq:nitsche}. 
Figure \ref{fig:nsym_p1_err} shows the the $H^1$ error $\norm{u - u_h}_{h,\star}$ and the $L^2$ error $\norm{ u - u_h}_{\Omega_h}$ for $\varepsilon=0,1$. 
We observe that the convergence rates are almost $O(h)$ for the $H^1$ error and $O(h^2)$ for the $L^2$ error. 
Thus, the optimal convergence rates actually take place and the estimates of Theorems \ref{thm1} and \ref{thm2} are confirmed.

\begin{figure}[bt]
  \centering
     \includegraphics[scale = 0.55]{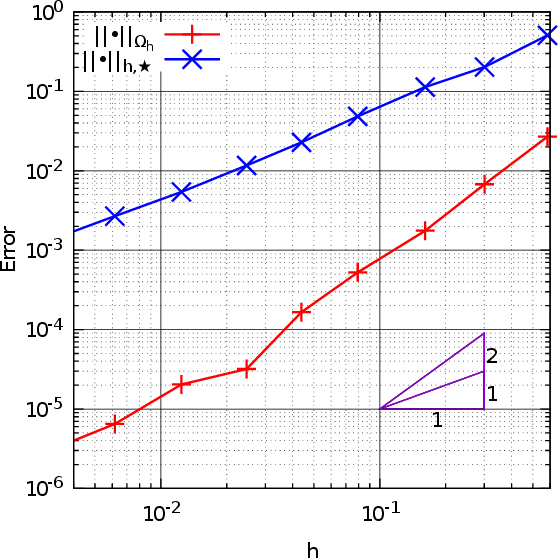}
  \caption{$H^1$ errors and $L^2$ errors for the exact solution $u=((x_1+1)^2+x_2^2)^{\frac{2}{3}}$}\label{fig:n_pow_err}
\end{figure}

Subsequently, we consider the exact solution $u(x_1,x_2) = ((x_1+1)^2+x_2^2)^{\frac{2}{3}}$ and the corresponding $f$, $g$, and $u_0$. Let $\varepsilon=1$.
In this case, $u \in H^2(\Omega)$ and $u \not\in W^{3,q}(\Omega)$ for any $q>n$.
That is, the assumption of Theorem \ref{thm2} does not satisfied.
Figure \ref{fig:n_pow_err} reports the $H^1$ error and the $L^2$ error for $\varepsilon=1$. We see from Figure \ref{fig:n_pow_err} that the convergence rate for the $H^1$ error is $O(h)$. However, the second-order convergence does not achieve for the $L^2$ error. Actually, we observe that the convergence rate for the $L^2$ error is $O(h^{2-\sigma})$ with some small $\sigma>0$.  
This result is consistent with Theorem \ref{thm2}. 
Therefore, we can conclude that the regularity condition $u\in W^{3,q}(\Omega)$ with $q>n$ is an essential requirement for deducing the optimal order convergence. 

\subsection*{Acknowledgements}
This work was supported by CREST (JPMJCR15D1) of JST, Japan, 
by Grant-in-Aid for Scientific Research B (15H03635) of JSPS, Japan, and 
by Grant-in-Aid for Scientific Research A (21H04431) of JSPS, Japan.


\end{document}